\RequirePackage{fix-cm}
\documentclass[smallextended]{svjour3}       
\smartqed  
\usepackage{graphicx}
\usepackage{amsfonts, amsmath, amscd}
\usepackage[psamsfonts]{amssymb}
\usepackage{amssymb}
\usepackage[usenames]{color}
\usepackage[all,cmtip]{xy}

\spnewtheorem{mytheorem}{Theorem}[section]{\bf}{\it}
\numberwithin{mytheorem}{section}
\spnewtheorem{mycorollary}[mytheorem]{Corollary}{\bf}{\it}
\spnewtheorem{mydefinition}[mytheorem]{Definition}{\bf}{}
\spnewtheorem{myproposition}[mytheorem]{Proposition}{\bf}{\it}
\spnewtheorem{myquestion}[mytheorem]{Question}{\bf}{\it}
\spnewtheorem{myproblem}[mytheorem]{Problem}{\bf}{\it}
\spnewtheorem{myremark}[mytheorem]{Remark}{\bf}{}
\spnewtheorem{myexample}[mytheorem]{Example}{\bf}{\it}
\spnewtheorem{myfact}[mytheorem]{Fact}{\bf}{\it}
\spnewtheorem{mylemma}[mytheorem]{Lemma}{\bf}{\it}

\newcommand{\CC}{C_k}

\newcommand{\w}{\omega}

\newcommand{\KK}{\mathcal{K}}

\newcommand{\IR}{\mathbb{R}}

\renewcommand{\phi}{\varphi}

\newcommand{\U}{\mathcal U}

\newcommand{\uhr}{\upharpoonright}

\begin{document}

\title{Topological properties of function spaces over ordinals\thanks{The second author was supported by the GACR project 15-34700L and RVO: 67985840. The third  author was supported by Generalitat Valenciana, Conselleria d'Educaci\'{o}, Cultura i Esport, Spain, Grant PROMETEO/2013/058 and by  the GA\v{C}R project 16-34860L and RVO: 67985840, and gratefully acknowledges also the financial support he received from the Kurt Goedel Research Center in Wien for his research visit in days 15.04-24.04 2016. The fourth author  would like to thank  the Austrian Science Fund FWF (Grant I 1209-N25) for generous support for this research. The collaboration of the second and the fourth authors was partially supported by the Czech Ministry of Education grant 7AMB15AT035 and RVO: 67985840.}
}

\titlerunning{Topological properties of function spaces over ordinals}        

\author{Saak Gabriyelyan \and  Jan Greb\'{\i}k \and  Jerzy K{\c{a}}kol \and Lyubomyr Zdomskyy
}

\authorrunning{S. Gabriyelyan \and  J. Greb\'{\i}k \and  J. K{\c{a}}kol \and L. Zdomskyy} 

\institute{S. Gabriyelyan  \at
              Department of Mathematics, Ben-Gurion University of the Negev, Beer-Sheva, P.O. 653, Israel \\
              \email{saak@math.bgu.ac.il}           
            \and
            J. Greb\'{\i}k \at
            Institute of Mathematics, Czech Academy of Sciences, Czech Republic \\
            \email{Greboshrabos@seznam.cz}
           \and
           J. K{\c{a}}kol \at
              A. Mickiewicz University $61-614$ Pozna{\'n}, Poland and Institute of Mathematics, Czech Academy of Sciences, Czech Republic\\
              \email{kakol@amu.edu.pl}
           \and
           L. Zdomskyy \at
           Kurt G\"odel Research Center for Mathematical Logic, University of Vienna, A-1090 Wien, Austria \\
            \email{lzdomsky@gmail.com}
}

\date{Received: date / Accepted: date}

\maketitle

\begin{abstract}
A topological space $X$ is said to be
an \emph{Ascoli space} if  any  compact subset $\KK$ of $\CC(Y)$ is evenly continuous.
This definition is motivated by the classical Ascoli theorem.
We study the $k_\IR$-property and the Ascoli property of $C_{p}(\kappa)$ and
 $\CC(\kappa)$ over ordinals  $\kappa$.  We prove that
$C_p(\kappa)$ is always an Ascoli space,   while  $C_p(\kappa)$ is a $k_\IR$-space iff
 the cofinality  of $\kappa$ is countable. In particular, this provides
the first  $C_{p}$-example of an Ascoli space which is not a  $k_\IR$-space, namely
$C_p(\w_1)$. We show that $\CC(\kappa)$ is Ascoli iff  $\mathrm{cf}(\kappa)$ is countable
 iff $\CC(\kappa)$ is metrizable.


\end{abstract}



\section{Introduction }


The study of topological properties   of  function spaces is
quite an  active area of research  attracting  specialists both from topology and functional analysis, see  for example \cite{Arhangel,BG,GKP,kak,Tkachuk-Book-1} and references therein. In the following diagram we select the most important compact type properties
 generalizing metrizability
\[
\xymatrix{
\mbox{metric} \ar@{=>}[r] & {\mbox{Fr\'{e}chet--}\atop\mbox{Urysohn}}  \ar@{=>}[r] & \mbox{$k$-space} \ar@{=>}[r] & \mbox{$k_\IR$-space} \ar@{=>}[r] & \mbox{Ascoli},}
\]
and note that none of these implications is reversible, see \cite{BG,Eng} (all relevant definitions are given in the next section).

For a Tychonoff  topological space $X$, we denote by $\CC(X)$ and $C_p(X)$ the space $C(X)$ of all continuous real-valued functions on $X$ endowed with the
compact-open topology and the  topology of pointwise convergence, respectively.

It is well-known that $C_p(X)$ is metrizable if and only if $X$ is countable.
Pytkeev, Gerlitz and Nagy (see \S 3 of \cite{Arhangel}) characterized  spaces $X$ for
which $C_p(X)$ is Fr\'{e}chet--Urysohn or a $k$-space (these properties coincide for
spaces of the form $C_p(X)$). The authors in \cite{GGKZ} obtained some sufficient conditions on $X$ for which the space $C_p(X)$ is an Ascoli space. Recall
 that $X$ is called an  {\em Ascoli space } if  any  compact subset $\KK$ of $\CC(X)$ is evenly continuous (or, equivalently, if the natural evaluation map $X\hookrightarrow \CC(\CC(X))$ is an embedding, see \cite{BG}).

Every linear order $<$ on a set $X$ generates a natural topology on $X$
whose subbase consists of sets of the form $\{z:z<x\}$ and $\{z:z>x\}$, where $x\in X$.
Spaces $X$ whose topology is generated by some linear order
 are called \emph{linearly ordered topological spaces.}
Ordinals with the topology generated by their natural wellorder
form is an interesting class of linearly ordered topological spaces, and
function spaces over them give a good source of (counter)examples in the
corresponding theory. For instance, the space $C_p(\w_1)$ is Lindel\"{o}f, see \cite{Tkachuk-Book-1}.
On the other hand,
Arhangel'skii showed in \cite{Arhangel-2001} that the space $C_p(\w_1+1)$
is not normal. In \cite{Gulko} Gul'ko proved  that there are no two distinct natural
number $n$ and $m$ for which the powers $C_p(\w_1)^n$ and $C_p(\w_1)^m$ are homeomorphic.
In \cite{Mor-Wul} Morris and Wulbert observed that $\CC(\w_1)$ is not barrelled.

In this short note we provide complete characterizations of those ordinals $\kappa$ for which $C_p(\kappa)$ and $\CC(\kappa)$ are $k_\IR$-spaces or Ascoli spaces. The following theorems are the main results of the paper.

\begin{mytheorem} \label{t:Ascoli-Cp-ordinal}
For every ordinal $\kappa$ the space $C_p(\kappa)$ is Ascoli.
\end{mytheorem}

Denote by $\mathrm{cf}(\kappa)$ the cofinality of an ordinal $\kappa$.
\begin{mytheorem} \label{t:k-R-Cp-ordinal}
For an ordinal $\kappa$, the space $C_p(\kappa)$ is a $k_\IR$-space if and only if
$\mathrm{cf}(\kappa)\leq\omega$ if and only if  $C_p(\kappa)$ is Fr\'{e}chet--Urysohn.
\end{mytheorem}

Theorems \ref{t:Ascoli-Cp-ordinal} and \ref{t:k-R-Cp-ordinal}  show
   that the space $C_p(\w_1)$ is Ascoli
but  not a $k_\IR$-space. This   answers   Question 6.8 in \cite{GKP} for spaces $C_p(X)$
in the affirmative, and complements \cite[3.3.E]{Eng}
asserting  that for uncountable discrete $X$ the space
$C_{p}(X)=\mathbb{R}^{X}$ is a $k_\IR$-space but   not a $k$-space.

\begin{mytheorem} \label{t:Ascoli-Ck-ordinal}
For an ordinal $\kappa$, the space $\CC(\kappa)$ is an Ascoli space if and only if $\mathrm{cf}(\kappa)\leq\omega$, so $\CC(\kappa)$ is complete and metrizable.
\end{mytheorem}

\section{Proofs}


Below we recall some topological  concepts  used in Theorems \ref{t:Ascoli-Cp-ordinal}
and \ref{t:Ascoli-Ck-ordinal}, for other notions we refer the reader to  the book \cite{Eng}.
A {\em $k$-cover} $\U$ of a topological space $X$ is a family of subsets
of $X$ such that every compact subset of $X$ is contained in some member of $\U$.

\begin{mydefinition} {\em
A topological space $X$ is

$\bullet$ {\em hemicompact} if there exists a countable $k$-cover of
$X$ consisting of compacts;

$\bullet$ {\em realcompact} if
it can be embedded  as a closed subset to $\mathbb R^\lambda$
for some cardinal $\lambda$;

$\bullet$  a {\em $k_\IR$-space} if  a real-valued function $f$ on $X$ is continuous if and only if its restriction $f|_K$ to any compact subset $K$ of $X$ is continuous;

$\bullet$ {\em scattered} if every non-empty subspace $A$ of $X$ has an isolated point in $A$.


}
\end{mydefinition}

Recall that an ordinal $\kappa$ is \emph{limit} if there is no $\alpha$ such that $\kappa=\alpha+1$,
otherwise $\kappa$ is called a \emph{successor} ordinal.
The {\em cofinality $\mathrm{cf}(\kappa)$} of a limit ordinal number $\kappa$ is the
 smallest ordinal $\alpha$ which is the order type of a cofinal subset of $\kappa$.
If $\kappa$ is a successor ordinal we set $\mathrm{cf}(\kappa)=1$.

 The following simple facts  should be well-known (for (i) see \cite[\S~5.11]{GiJ},
the other ones are straightforward). 
\begin{mylemma} \label{l:order-pseudocom}
\begin{enumerate}
\item[{\rm (i)}] $\kappa$ is compact if and only if it is a successor;
\item[{\rm (ii)}] $\kappa$ is hemicompact non-countably compact if and only if $\mathrm{cf}(\kappa)=\w$;
\item[{\rm (iii)}] $\kappa$ is countably compact non-compact if and only if $\mathrm{cf}(\kappa)>\w$.
\end{enumerate}
\end{mylemma}

For the convenience of the reader we recall also the following two results.
\begin{myproposition}[\cite{GKP}] \label{p:Ascoli-sufficient}
Assume   $X$ admits a  family $\U =\{ U_i : i\in I\}$ of open subsets of $X$, a subset $A=\{ a_i : i\in I\} \subset X$ and a point $z\in X$ such that: {\rm (i)} $a_i\in U_i$ for every $i\in I$, {\rm (ii)} $\big|\{ i\in I: C\cap U_i\not=\emptyset \}\big| <\infty$  for each compact subset $C$ of $X$, and {\rm (iii)} $z$ is a cluster point of $A$. Then $X$ is not an Ascoli space.
\end{myproposition}


A family $\{ A_i\}_{i\in I}$ of subsets of a set $X$ is said to be {\em  point-finite} if the set $\{i\in I: x\in A_i\}$ is finite for every $x\in X$.
A family $\{ A_i\}_{i\in I}$ of subsets of a topological space $X$ is called {\em strongly point-finite} if for every $i\in I$, there exists an open set
$U_i$ of $X$ such that $A_i\subseteq U_i$ and $\{ U_i\}_{i\in I}$ is point-finite.  Following Sakai \cite{Sak2}, a topological space $X$ is said to have the
{\em property $(\kappa)$} if every pairwise disjoint sequence of finite subsets  of $X$ has a strongly point-finite subsequence.
The following result is proved in \cite{GGKZ}.

\begin{mytheorem} \label{t:Cp-Ascoli-k-FU}
If $C_p(X)$ is Ascoli,  then $X$ has the property $(\kappa)$.
\end{mytheorem}

It is well-known that ordinals are locally compact and scattered (for the last property we note that the smallest element of a subset $A$ of $X$ is isolated in $A$). The following proposition is of independent interest, it generalizes Corollary 1.5 of \cite{GGKZ} and immediately implies Theorem \ref{t:Ascoli-Cp-ordinal}.
\begin{myproposition} \label{p:locally-compact-Ascoly}
Let $X$ be a locally compact space. Then $C_{p}(X)$ is Asoli if and only if $X$ is scattered.
\end{myproposition}

\begin{proof}
The ``only if'' part follows from Theorem~\ref{t:Cp-Ascoli-k-FU}
combined with the fact that the property $(\kappa)$ is preserved by subspaces, along with the fact that every compact space with the property $(\kappa)$ is scattered, see \cite[Theorem~3.2]{Sak2}. For the ``if'' direction consider the one-point compactification $X^*=X\cup\{x_\infty\}$   of $X$ and note that it is scattered. Therefore  $C_p(X^*)$ is Frechet-Urysohn by \cite[II.7.16]{Arhangel}, and hence so is its subspace  $Z$ consisting of those continuous $f:X^*\to\IR$ such that $f(x_\infty)=0$. Now it is easy to see that $Z\uhr X=\{f\uhr X:f\in Z\}\subset C_p(X)$ is homeomorphic to $Z$ and is dense in $C_p(X)$. Therefore $C_p(X)$ has a
 dense Frechet--Urysohn subspace, and hence every function $f$ belongs to a dense Ascoli subspace $f+Z$ of $C_p(X)$. Thus $C_p(X)$ is Ascoli by Proposition 5.10 of \cite{BG}. $\Box$
\end{proof}

Below we prove Theorems \ref{t:k-R-Cp-ordinal} and \ref{t:Ascoli-Ck-ordinal}.

{\em Proof of Theorem \ref{t:k-R-Cp-ordinal}.}
Let $C_p(\kappa)$ be a $k_\IR$-space and suppose towards
a contradiction that $\mathrm{cf}(\kappa)>\omega$.
Then $\kappa$ is countably compact by Lemma \ref{l:order-pseudocom}. Hence the space $C(\kappa)$
endowed with the sup-norm is a Banach space. Therefore the space $C_p(\kappa)$ admits a stronger normed topology and is angelic by \cite[Proposition 9.6]{kak}.
 Since every compact subset of $C_p(\kappa)$ is Fr\'{e}chet--Urysohn, every sequentially
continuous function on $C_p(\kappa)$ is continuous. In particular, every sequentially continuous
 functional on $C_p(\kappa)$ is continuous. So $\kappa$ is a realcompact space by
Theorem 1.1 of \cite{Wilansky}. Being realcompact and pseudocompact the space $\kappa$
is compact by \cite[3.11.1]{Eng}.
Hence $\kappa$ is a successor ordinal, a contradiction. Thus $\mathrm{cf}(\kappa)\leq\w$.

Conversely, let $\mathrm{cf}(\kappa)\leq\w$. Then $\kappa$ is hemicompact by
Lemma \ref{l:order-pseudocom}. Thus $C_p(\kappa)$ is Fr\'{e}chet--Urysohn
by \cite[II.7.16]{Arhangel}.
$\Box$

Since each ordinal $\alpha$ is  the set of all smaller ordinals,
in the following proof
we adopt the following perhaps standard notation:
For a function $f$ whose  domain is an ordinal $\kappa$
and $\alpha\in\kappa$ we denote by $f(\alpha)$ the value of $f$ at $\alpha$,
and by $f[\alpha]$ the set $\{f(\beta):\beta<\alpha\}$.
\medskip

{\em Proof of Theorem \ref{t:Ascoli-Ck-ordinal}.}
Suppose for a contradiction that $\mathrm{cf}(\kappa)>\omega$.
We shall use Proposition \ref{p:Ascoli-sufficient}
and show that $C_k(\kappa)$ is not Ascoli.
For every $\alpha<\kappa$ we define $f_{\alpha}:\kappa \to [0,1]$
by $f_\alpha\big[\alpha+1\big]=\{0\}$ and $f_\alpha\big[\kappa\setminus(\alpha+1)\big]=\{1\}$, and set
\[
U_{\alpha}:=\{f\in \CC(\kappa): \; f(\alpha)<1/4, \ f(\alpha+1)>3/4\}.
\]
To prove that $\CC(\kappa)$ is not Ascoli it is enough to verify
 the assumptions of Proposition \ref{p:Ascoli-sufficient} for
$\{f_{\alpha}\}_{\alpha<\kappa}$, $\{U_{\alpha}\}_{\alpha<\kappa}$ and $0\in \CC(\kappa)$.
 Clearly, (i) and (iii) hold true. Let us check (ii).
 Take  any compact $C\subseteq \CC(\kappa)$ and assume, contrary to our claim,
that there are infinitely many $\alpha<\kappa$ such that $C\cap U_\alpha\neq \emptyset$.
Then there exists a strictly increasing sequence  $\{\alpha_n\}_{n<\omega}$
 such that $C\cap U_{\alpha_n}\not=\emptyset$.
Let $\alpha=\lim \alpha_n$. As $\mathrm{cf}(\kappa)>\omega$ we have
$\alpha<\kappa$. By the Ascoli theorem used for $\alpha+1\subset \kappa$ and
$1/ 2$ we can find  a basic neighborhood
$O_{\alpha}$ of $\alpha$ such that $|h(x)-h(y)|<1/4$ for all $x,y\in O_{\alpha}$ and
$h\in C$. Take $n$ such that $\alpha_n\in O_\alpha$
and fix $h\in C\cap U_{\alpha_n}$. Then
\begin{eqnarray*}
\begin{aligned}
\frac{1}{4} & > |h(\alpha_n+1)-h(\alpha_n)| \\
& \geq |f_{\alpha_n}(\alpha_n+1)-f_{\alpha_n}(\alpha_n)|-
|f_{\alpha_n}(\alpha_n+1)-h(\alpha_n+1)| -|h(\alpha_n)-f_{\alpha_n}(\alpha_n)| \\
& >1-\frac{1}{4}-\frac{1}{4}= \frac{1}{2},
\end{aligned}
\end{eqnarray*}
 which is a contradiction. Thus $\mathrm{cf}(\kappa)\leq\omega$.

Conversely, if $\mathrm{cf}(\kappa)\leq\omega$, then $\kappa$ is a
hemicompact locally compact space by Lemma \ref{l:order-pseudocom}. Hence $\CC(\kappa)$ is complete metrizable by Corollary 5.2.2 of \cite{mcoy}.
$\Box$

\bibliographystyle{amsplain}

\end{document}